\documentclass[final,5p,times,twocolumn]{elsarticle}

\usepackage{preamble}

\usetikzlibrary{decorations.pathreplacing, positioning}
\journal{European Journal of Control}

\begin{document}

\begin{frontmatter}

\title{Scaled Relative Graphs for Nonmonotone Operators with Applications in Circuit Theory}
\tnotetext[labelfund]{This work was supported by the Research Foundation Flanders (FWO) PhD grant 1183822N and research projects G081222N, G033822N, and G0A0920N;
            Research Council KUL grant C14/24/103; and the European Research Council under the European Union's Horizon 2020 research and innovation program / ERC Advanced Grant: SpikyControl (no.\ 101054323).}

\author{Jan Quan\corref{cor1}}
\ead{jan.quan@kuleuven.be}
\author{Brecht Evens}
\ead{brecht.evens@kuleuven.be}
\author{Rodolphe Sepulchre}
\ead{rodolphe.sepulchre@kuleuven.be}
\author{Panagiotis Patrinos}
\ead{panos.patrinos@kuleuven.be}
\cortext[cor1]{Corresponding author.}

\affiliation{organization={KU Leuven, Department of Electrical Engineering ESAT-STADIUS},%
            addressline={Kasteelpark Arenberg 10, box 2446}, 
            city={Leuven},
            postcode={3001}, 
            country={Belgium}}

\begin{abstract}
The scaled relative graph (SRG) is a powerful graphical tool for analyzing the properties of operators, by mapping their graph onto the complex plane.
In this work, we study the SRG of two classes of nonmonotone operators, namely the general class of semimonotone operators and a class of angle-bounded operators. 
In particular, we provide an analytical description of the SRG of these classes and show that membership of an operator to these classes can be verified through geometric containment of its SRG.
To illustrate the importance of these results, we provide several examples in the context of electrical circuits.
Most notably, we show that the Ebers--Moll transistor belongs to the class of angle-bounded operators and use this result to compute the response of a common-emitter amplifier using Chambolle--Pock, despite the underlying nonsmoothness and multi-valuedness, leveraging recent convergence results for this algorithm in the nonmonotone setting. 

\end{abstract}


\end{frontmatter}

\section{Introduction}
	\label{sec:introduction}
	
Recently, the scaled relative graph (SRG) has emerged as a powerful tool for analyzing individual operators and their interconnections. 
Originally introduced in \cite{ryu2022scaled}, the SRG can be interpreted as a generalization of the classical Nyquist diagram to arbitrary nonlinear operators.
By mapping an operator's graph onto the complex plane, the SRG provides insights into its incremental gain and phase properties, with interconnections represented as graphical manipulations of SRGs.
For instance, this approach unifies and extends classical results like the Nyquist criterion and the incremental passivity theorem~\cite{chaffey2023graphical}, and can be used as a formal framework for constructing geometric proofs of convergence for contractive and nonexpansive fixed-point iterations, with similar 2D visualizations appearing earlier in \cite{eckstein1992DouglasRachford,giselsson2016Linear}.
Central to this approach is the concept of SRG-full operator classes, where membership of an operator to such classes directly corresponds to the geometric containment of its SRG.
SRG-full classes include many common operator classes defined through inequalities, such as (hypo)monotone, co(hypo)monotone, Lipschitz, and averaged operators \cite{ryu2022scaled}.

In many applications, operators belonging to these classes emerge quite naturally.
For instance, in circuit theory many commonly used circuit elements such as linear time-invariant resistors, capacitors, inductors, transformers and gyrators are known to be maximally monotone \cite{chaffey2023circuit,chaffey2024monotone}.
Based on this observation, \cite{chaffey2023circuit} demonstrated that the behavior of a monotone circuit can be modeled as the zero of a monotone + skew inclusion problem, where the monotone component represents the device equations and the skew component arises from the circuit topology by Tellegen's theorem \cite{tellegen1952general}. 
This type of inclusion problem can be solved in an efficient manner using the Chambolle--Pock iteration \cite{chambolle2011firstorder}.
Compared to traditional methods based on ordinary differential equations, this splitting method offers greater scalability, robustness to parameter variations, and is able to deal with nonsmooth and multi-valued elements.

However, besides these monotone circuit elements, there are also several commonly used circuit elements which exhibit complicated \emph{nonmonotone} behavior, not properly captured by any previously mentioned SRG-full operator class.
Notable examples include the tunnel diode, transistor devices, nonlinear capacitors/inductors and memristive elements.
For instance, using the SRG it was observed in \cite{chaffey2024monotone} that the potassium conductance in the Hodgkin--Huxley membrane model \cite{hodgkin1952components} is hypomonotone, although it also exhibits stronger properties which this class does not adequately capture.

In this work, we will address this issue by considering two classes of nonmonotone operators which are more suitable for accurately capturing the nonmonotone behavior of these circuit elements, namely \emph{semimonotone operators} and \emph{angle-bounded operators}.
The class of semimonotone operators was first introduced in \cite{evens2023convergence} and can be used to derive sufficient conditions for several splitting methods in the nonmonotone setting, including Chambolle--Pock \cite{evens2023convergenceCP}.
This class was recently studied in~\cite{huijzer2025modelling} for networks consisting only of memristors.
The class of angle-bounded operators is inspired by the singular angle introduced in \cite{chen_singular_2021}.
The singular angle is the phase counterpart of the $L_2$-gain, capturing the amount of rotation induced by a system.

Our main contributions are as follows:
\begin{enumerate}[label=(\roman*)]
	\item We derive an analytical expression for the SRGs of semimonotone and angle-bounded operators, and show that both are SRG-full. Additionally, we establish a connection between these two classes of operators through the SRG.
	\item 
	We show that Ebers--Moll transistors are angle-bounded under standard assumptions, and as a result also semimonotone.
	\item 
	We consider common-emitter amplifiers involving transistors and tunnel diodes.
	Despite the nonsmoothness, nonmonotonicity, and potential multi-valuedness of these elements, we show that their response can be computed efficiently using the Chambolle--Pock algorithm, building on recent convergence results for semimonotone operators \cite{evens2023convergenceCP}.
	Due to the aforementioned difficulties, standard Newton methods are not applicable in this setting, highlighting the merits of our results.
\end{enumerate}

\subsection{Notation}
	We denote the set of complex and extended-complex numbers by \(\mathbb{C}\) and \(\overline{\mathbb{C}} \coloneqq \mathbb{C} \cup \{\infty\}\) respectively. The graph of a set-valued mapping $A : \mathcal{H} \rightrightarrows \mathcal{H}$ on a Hilbert space $\mathcal{H}$ is defined as $\graph A \coloneqq \{(x,y) \in \mathcal{H}\times\mathcal{H} \mid y \in A(x)\}$, and 
$A$ is said to be outer semicontinuous if its graph is a closed subset of $\mathcal{H}\times\mathcal{H}$.
An operator class $\mathcal{A}$ is a set of operators on Hilbert spaces.
We denote the identity operator on a suitable space by $\id$.
For scalars and sets, Minkowksi-type operations are to be understood, i.e.,
$A + B \coloneqq \{a + b \mid a \in A, b \in B\}$
and
$\alpha A \coloneqq \{\alpha a \mid a \in A\}$.
An open disk with center $c \in \mathbb{C}$ and radius $r > 0$ is defined as $D(c,r) \coloneqq \{z \in \mathbb{C} \mid \|z - c \| < r\}$. The $2$-argument arctangent function $\atantwo(y,x)$ denotes the phase of the complex number $x+iy$, confined to $(-\pi, \pi]$.

\section{Scaled relative graphs}\label{sec:srg}

First, we introduce the concept of scaled relative graphs, which map the incremental properties of an operator $A : \mathcal{H}\rightrightarrows\mathcal{H}$ to a subset of the extended complex plane $\overline{\mathbb{C}}$.

Consider a pair $(x,u),(y,v) \in \graph A$ and define the corresponding complex conjugate pair
\[
    z_{\pm}(x-y, u-v) \coloneqq \frac{\|u-v\|}{\|x-y\|}\exp(\pm i \angle (x-y, u-v)),
\]
where 
\[
    \angle(x-y, u-v)
        \coloneqq
    \begin{cases}
        \arccos \left(\frac{\langle x-y, u-v\rangle}{\|x-y\|\|u-v\|}\right) & \textnormal{if $x\neq y$ and $u\neq v$}, \\
        0 & \textnormal{otherwise}.
    \end{cases}    
\]
By considering all different pairs, the scaled relative graph of operators and operator classes can be constructed as follows.
\begin{definition}\label{def:SRG}
    The SRG of an operator $A$ is defined as
    \begin{multline*}
        G(A) \coloneqq \left\{z_{\pm}(x-y, u-v) \mathrel{\big|\; u\in A(x), v\in A(y), x\neq y}\right\} \\
        \Big(\cup \{\infty\} \textnormal{ if $A$ is multi-valued}\Big)
    \end{multline*}
    and the SRG of an operator class $\operatorclass{A}$ is defined as
    \[
        G(\operatorclass{A})
            \coloneqq
        \bigcup_{A \in \operatorclass{A}} G(A).
    \]
\end{definition}
By construction, the SRG provides a visualization of the incremental gain and phase of all input-output pairs of an operator (or an operator class).
Note that such a visualization is similar to the Nyquist diagram of a linear time-invariant (LTI) transfer function.
In fact, the SRG of a stable LTI transfer function is the convex hull of its Nyquist diagram, under the
Beltrami-Klein mapping \cite[Thm.\ 4]{chaffey2023graphical}.
For more details on this connection, we refer the interested reader to \cite{chaffey2021scaled,pates2021scaled}.

The main focus of this paper is on SRG-full operator classes, defined as follows.
\begin{definition} \label{def:srgfull}
    An operator class $\operatorclass{A}$ is SRG-full if
    \[
        A \in \operatorclass{A} \Leftrightarrow G(A) \subseteq G(\operatorclass{A}).
    \]
\end{definition}
In essence, a class is SRG-full if membership of an operator to that class is equivalent to the geometric containment of its SRG, enabling class membership to be verified graphically.
In particular, we will focus on SRG-full classes defined by a nonnegatively homogeneous function $h$, i.e., a function for which $h(\eta a,\eta b,\eta c) = \eta h(a,b,c)$ for all $\eta \geq 0$.
\begin{proposition} \label{prop:srgfull} (\hspace{1sp}\cite[Thm.\ 2]{ryu2022scaled})
    An operator class $\operatorclass{A}$ is SRG-full if there is a nonnegatively homogeneous function $h : \mathbb{R}^3 \to \mathbb{R}$ such that 
    \begin{multline*}
        A \in \operatorclass{A} \Leftrightarrow  h(\|u-v\|^2, \|x-y\|^2, \langle x-y, u-v\rangle) \leq 0, \\
        \forall (x,u), (y,v) \in \graph A.
    \end{multline*}
\end{proposition}
Notable examples of such operator classes include the classes of (hypo)monotone, co(hypo)monotone, Lipschitz, and averaged operators \cite{ryu2022scaled}.

Before proceeding, we recall some important calculus rules for SRGs which will be used throughout the paper.
\begin{proposition} \label{prop:srgcalculus} (\hspace{1sp}\cite[Thms.\ 4\&5]{ryu2022scaled})
    Let $\operatorclass{A}$ be an operator class, and let $\alpha \in \mathbb{R} \setminus \{0\}$. Then, the following hold.
    \begin{enumerate}
        \item \label{prop:srgcalculus:1}
        \(        
            G(\alpha \operatorclass{A}) = G(\operatorclass{A} \alpha) = \alpha G(\operatorclass{A})
        \)
        \item \label{prop:srgcalculus:2}
        \(        
            G(\textnormal{id} + \operatorclass{A}) = 1 + G(\operatorclass{A})
        \)
        \item \label{prop:srginv}
        \(
            G(\operatorclass{A}^{-1}) = (G(\operatorclass{A}))^{-1} = \{\frac1r e^{i\varphi} \mid r e^{i\varphi} \in G(\operatorclass{A})\}
        \)
    \end{enumerate}
\end{proposition}

\begin{figure*}[hbpt]
    \centering
    \begin{tabular}{@{}c@{\;\;}c@{\;\;}c@{\;\;}c@{\;\;}c@{}}
        \begin{tikzpicture}[baseline=(current bounding box.south)]
\def\signrho{0} %
\def\rhonum{0}
\def\signsigma{0} 
\def\sigmanum{0}

\input{figures/draw_semimonotone.tex}
\end{tikzpicture}&
        \begin{tikzpicture}[baseline=(current bounding box.south)]
\def\signrho{0} %
\def\rhonum{0}
\def\signsigma{-1} 
\def\sigmanum{-0.5}

\input{figures/draw_semimonotone.tex}
\end{tikzpicture}&
        \begin{tikzpicture}[baseline=(current bounding box.south)]
\def\signrho{0} %
\def\rhonum{0}
\def\signsigma{1} 
\def\sigmanum{0.5}

\input{figures/draw_semimonotone.tex}
\end{tikzpicture}&
        \begin{tikzpicture}[baseline=(current bounding box.south)]
\def\signrho{-1} %
\def\rhonum{-0.75}
\def\signsigma{0} 
\def\sigmanum{0}

\input{figures/draw_semimonotone.tex}
\end{tikzpicture}&
        \begin{tikzpicture}[baseline=(current bounding box.south)]
\def\signrho{1} %
\def\rhonum{0.75}
\def\signsigma{0} 
\def\sigmanum{0}

\input{figures/draw_semimonotone.tex}
\end{tikzpicture}\\
        $G(\operatorclass{M})$ &
        $G(\operatorclass{M}_{-\mu})$ & $G(\operatorclass{M}_{\mu})$ & $G(\operatorclass{C}_{-\rho})$ & $G(\operatorclass{C}_{\rho})$
    \end{tabular}
    \caption[
        Scaled relative graphs for the classes of hypomonotone, strongly monotone, cohypomonotone and cocoercive operators.
    ]{
        Scaled relative graphs for the classes of monotone ($\operatorclass{M}$),
        hypomonotone ($\operatorclass{M}_{-\mu}$), strongly monotone ($\operatorclass{M}_{\mu}$), cohypomonotone ($\operatorclass{C}_{-\rho}$) and cocoercive operators ($\operatorclass{C}_{\rho}$), where $\mu > 0$ and $\rho > 0$.
        The shaded regions indicate portions of the extended complex plane included in the corresponding SRG.
    }
    \label{fig:semi:srg:simple}
\end{figure*}

As an example of these calculus rules, consider the SRG of the class of monotone operators, given by
\[
    G(\operatorclass{M}) = \{z \in \mathbb{C} \mid \realpart z \geq 0\} \cup \{\infty\}.
\]
By \Cref{prop:srgcalculus}, it follows immediately that the SRGs of the class of $\mu$-monotone operators $\operatorclass{M}_{\mu} = \mu \textnormal{id} + \operatorclass{M}$ and $\rho$-comonotone operators $\operatorclass{C}_{\rho} = \operatorclass{M}_{\rho}^{-1}$ are given by
\begin{align}
    G(\operatorclass{M}_\mu) 
        {}={}&
    \{z \in \mathbb{C} \mid \realpart z \geq \mu\}\cup \{\infty\},\\
    G(\operatorclass{C}_{\rho})
        {}={}& 
    \begin{cases}
        \overline{\mathbb{C}} \setminus D\left(\frac{1}{2\rho}, \frac{1}{2|\rho|}\right) & \textnormal{if $\rho < 0$}, \\
        \closure D\left(\frac{1}{2\rho}, \frac{1}{2|\rho|}\right) & \textnormal{if $\rho > 0$},
    \end{cases}
    \label{eq:SRG:comon}
\end{align}
where $\closure$ denotes the closure of a set.
A visualization of these SRGs is provided in \Cref{fig:semi:srg:simple}.

\section{SRGs of nonmonotone operators}\label{sec:srg-nonmon}
	In this section, we examine the SRG of two recently introduced classes of nonmonotone operators, namely semimonotone and angle-bounded operators.
We start by defining the semimonotone class, first introduced in \cite[Def.\ 4.1]{evens2023convergence}.
\begin{definition} 
    Let $\mu, \rho \in \mathbb{R}$. A set-valued operator $A : \mathcal{H} \rightrightarrows \mathcal{H}$ is ($\mu, \rho$)-semimonotone if
    \begin{align*}
        \langle x-y, u-v\rangle \geq \mu \|x-y\|^2 + \rho\|u-v\|^2, \quad \forall (x, u), (y, v) \in \graph A.
    \end{align*}
    It is maximally $(\mu,\rho)$-semimonotone if its graph is not strictly contained in the graph of any other $(\mu,\rho)$-semimonotone operator. The set of $(\mu, \rho)$-semimonotone operators will be denoted by $\operatorclass{S}_{\mu,\rho}$.
\end{definition}
Since this class is defined through the nonnegatively homogeneous function $h : (a,b,c) \mapsto \rho a + \mu b - c$, it is SRG-full by \Cref{prop:srgfull}.
\begin{proposition} \label{prop:semisrgfull}
    $\operatorclass{S}_{\mu,\rho}$ is SRG-full for any $\mu,\rho \in \mathbb{R}$.
\end{proposition}
Note that this operator class generalizes many well-known operator classes in the literature, including monotone and comonotone operators ($\operatorclass{S}_{0,0} = \operatorclass{M}$, $\operatorclass{S}_{\mu,0} = \operatorclass{M}_\mu$ and $\operatorclass{S}_{0,\rho} = \operatorclass{C}_{\rho}$). Equivalences with other operator classes are detailed in \cite[Rem.\ 4.2]{evens2023convergence}.

Now, we recall a calculus rule for the sum of a semimonotone operator with identity from \cite[Prop.\ 4.8]{evens2023convergence}, which is then used alongside \Cref{prop:srgcalculus} to obtain the analytical expression for the SRG of a semimonotone operator. 

\begin{proposition} (\hspace{1sp}\cite[Prop.\ 4.8]{evens2023convergence}) \label{prop:addidentitysm}
    Let $A : \mathcal{H} \rightrightarrows \mathcal{H}$ be (maximally) $(\mu, \rho)$-semimonotone and consider $B = A + \alpha \id$ with $\alpha\in\mathbb{R}$. If $1+2\rho \alpha > 0$, then $B$ is (maximally) $\left(\frac{\mu + \alpha(1+\rho\alpha)}{1+2\rho\alpha}, \frac{\rho}{1+2\rho\alpha}\right)$-semimonotone.
\end{proposition}
\begin{proposition}
    \label{prop:srgsemi}
    Let $\mu, \rho \in \mathbb{R}\setminus\{0\}$ such that $\mu\rho < \frac{1}{4}$. Let $c=\frac{1}{2\rho}$ and $r= \frac{\sqrt{1-4\mu\rho}}{2|\rho|}$. Then the SRG of the operator class of $(\mu,\rho)$-semimonotone operators is given by 
    \[
        G(\operatorclass{S}_{\mu,\rho}) = \begin{cases}
            \overline{\mathbb{C}} \setminus D(c, r) & \textnormal{if $\rho < 0$}, \\
            \closure D(c, r) & \textnormal{if $\rho > 0$}.
        \end{cases}
    \]
\end{proposition}
\begin{proof}
    Let $\alpha = \frac{-1 + \sqrt{1-4\mu\rho}}{2\rho}$ and note that $1+2\rho\alpha > 0$. Consequently, it follows from \Cref{prop:addidentitysm} that $\operatorclass{S}_{\mu,\rho} + \alpha \id = \operatorclass{S}_{0, \beta} = \operatorclass{C}_{\beta}$ where $\beta = \frac{\rho}{1+2\rho\alpha}$.
    The result then immediately follows from \eqref{eq:SRG:comon} and \Cref{prop:srgcalculus:1,prop:srgcalculus:2}, using the fact that $\frac{1+2\rho\alpha}{2\rho} - \alpha = \frac{1}{2\rho}$ and $\frac{1+2\rho\alpha}{2|\rho|} = \frac{\sqrt{1-4\mu\rho}}{2|\rho|}$.
\end{proof}
There are four qualitatively different SRGs for semimonotone operators, depending on the signs of $\mu$ and $\rho$. These are shown in \Cref{fig:semi:srg:semi}.

\begin{figure*}[hbpt]
    \centering
    \begin{tabular}{@{}c@{\;\;}c@{\;\;}c@{\;\;}c@{\;\;}c@{}}
        \begin{tikzpicture}[baseline=(current bounding box.south)]
    \def\signrho{-1} %
    \def\rhonum{-0.8}
    \def\signsigma{-1}
    \def\sigmanum{-0.1}

    \input{figures/draw_semimonotone.tex}
    \end{tikzpicture}&
        \begin{tikzpicture}[baseline=(current bounding box.south)]
    \def\signrho{-1} %
    \def\rhonum{-0.8}
    \def\signsigma{1}
    \def\sigmanum{0.1}

    \input{figures/draw_semimonotone.tex}
    \end{tikzpicture}&
        \begin{tikzpicture}[baseline=(current bounding box.south)]
    \def\signrho{1} %
    \def\rhonum{0.8}
    \def\signsigma{1}
    \def\sigmanum{0.1}

    \input{figures/draw_semimonotone.tex}
    \end{tikzpicture}&
        \begin{tikzpicture}[baseline=(current bounding box.south)]
    \def\signrho{1} %
    \def\rhonum{0.8}
    \def\signsigma{-1}
    \def\sigmanum{-0.1}

    \input{figures/draw_semimonotone.tex}
    \end{tikzpicture}\\
        $G(\operatorclass{S}_{-\mu,-\rho})$ & $G(\operatorclass{S}_{\mu,-\rho})$ & $G(\operatorclass{S}_{\mu,\rho})$ & $G(\operatorclass{S}_{-\mu,\rho})$
    \end{tabular}
    \caption[
        Four qualitatively different scaled relative graphs for classes of semimonotone operators.
    ]{
        Four qualitatively different scaled relative graphs for classes of semimonotone operators, where $\mu > 0$ and $\rho > 0$.
        The circles in the first and third graphs have radius $\nicefrac{\sqrt{1-4\mu\rho}}{2|\rho|}$, while those in the second and fourth have radius $\nicefrac{\sqrt{1+4\mu\rho}}{2|\rho|}$, as established in \Cref{prop:srgsemi}.
        }
    \label{fig:semi:srg:semi}
\end{figure*}

\begin{remark}
    If $\mu < 0$, $\rho < 0$ and $\mu\rho \geq \frac{1}{4}$, then all operators are $(\mu,\rho)$-semimonotone.
    Similarly, if $\mu > 0$, $\rho > 0$ and $\mu\rho > \frac{1}{4}$, then there exist no operators that are $(\mu,\rho)$-semimonotone \cite[Prop.\ 4.3]{evens2023convergence}. 
\end{remark}
Having established the SRG of the class of semimonotone operators, we continue with a class of angle-bounded operators, defined as follows.
\begin{definition} \label{def:anglebounded}
    Let $\theta \in [0, \pi]$.
    A set-valued operator $A : \mathcal{H} \rightrightarrows \mathcal{H}$ is $\theta$-angle-bounded if
    \[
        \angle (x-y, u-v) \leq \theta, \quad \forall (x,u), (y,v) \in \graph A,
    \]
    or equivalently if
    \begin{align}
        \langle x-y, u-v \rangle \geq \cos(\theta) \| x-y \| \| u-v \|.
        \label{eq:def:anglebounded:alt}
    \end{align}
    The set of all $\theta$-angle-bounded operators is denoted by 
    \(
        \operatorclass{B}_{\theta}.
    \)
\end{definition}
This notion is closely related to the singular angle defined in~\cite{chen_singular_2021}.
By construction, the SRG of the class of $\theta$-angle-bounded operators is given by
\[
    G(\operatorclass{B}_{\theta}) = \{re^{i\tilde{\theta}} \mid r \in \mathbb{R}_+, \tilde{\theta} \in [-\theta, \theta] \}\cup \{\infty\}.
\]
This SRG is visualized in \Cref{fig:anglebounded}. Note that this operator class is also SRG-full as it is defined through the nonnegatively homogeneous function $h(a,b,c) = \cos(\theta)\sqrt{a}\sqrt{b} - c$.

In what follows, we show that an angle-bounded operator, as well as the sum of an angle-bounded operator and an identity operator, is semimonotone. The proof follows via a geometric argument, leveraging the SRG-full property of semimonotone operators.
Our interest in this result stems from the fact that, in certain applications, verifying angle-boundedness may be more straightforward than verifying semimonotonicity directly.
For example, we will later use this result in \Cref{cor:transistor:semi} to show that the Ebers--Moll transistor is semimonotone under classical assumptions.

\begin{proposition} \label{prop:anglesem}
    Let $\rho<0$, $\alpha\geq0$, $\theta\in (\frac{\pi}{2}, \pi)$ and define
    \begin{align}
        \mu = \frac{1-(1-2\alpha\rho)^2\sin^2(\theta)}{4\rho}.
        \label{SRG:eq:optimal-mu}
    \end{align}
    Then, $\operatorclass{B}_{\theta} + \alpha \id \subseteq \operatorclass{S}_{\mu,\rho}$.
\end{proposition}
\begin{proof}
    By \Cref{prop:semisrgfull,prop:srgsemi}, it holds that $\operatorclass{B}_{\theta} + \alpha \id \subseteq \operatorclass{S}_{\mu,\rho}$ if and only if $G(\operatorclass{B}_{\theta} + \alpha \id)$ does not contain the open disk with center $\nicefrac{1}{2\rho}$ and radius
    \(
        r = \tfrac{\sqrt{1-4\mu\rho}}{2|\rho|} 
    \)
    as visualized in \Cref{fig:semiangleproof}. By the law of sines, this disk is not contained if
    \[
        r
        \leq
        \sin(\theta)\left(\frac{1}{2|\rho|} + \alpha\right),
    \]
    from which the claimed result follows by simple algebraic manipulation.
\end{proof}

Note that when invoking \Cref{prop:anglesem} for a given operator class $\operatorclass{B}_\theta + \alpha \id$, the choice of $\rho < 0$ provides a degree of freedom in how we express its semimonotonicity. This freedom becomes apparent through the geometric argument in the proof and is visualized in \Cref{fig:semiangle:non-unique}. However, when $\alpha > 0$, we may restrict our attention to the class of comonotone operators by setting $\mu = 0$ in \eqref{SRG:eq:optimal-mu}, leading to the following corollary.

\begin{cor}\label{cor:anglesem:identity}
    Let $\theta \in (\frac{\pi}{2}, \pi)$ and $\alpha > 0$. Then, $\operatorclass{B}_{\theta} + \alpha \id \subseteq \operatorclass{C}_{\rho}$, where $\rho = \frac{1}{2\alpha}(1 - \frac{1}{\sin \theta}) < 0$.
\end{cor}

\begin{figure}
    \centering
    \hspace{-0.2cm}
    \subfloat[]{%
        \begin{tikzpicture}[baseline=(current bounding box.south)]
    \pgfkeys{/pgf/number format/.cd,frac,frac whole=false}
    \pgfmathsetmacro\radius{0.95};
    \pgfmathsetmacro\angle{120};

    \pgfmathsetmacro\boxlength{1.4*\radius}
    \pgfmathsetmacro\gridlength{1.4*\radius}
    \pgfmathsetmacro\legendlength{1.25*\radius}

    \useasboundingbox (-\boxlength,-\boxlength) rectangle (\boxlength,\boxlength);
    \clip (-\boxlength,-\boxlength) rectangle (\boxlength,\boxlength);

    \draw[fill=firstcol, draw=white, very thick, fill opacity = 0.2, draw opacity = 0] (0,0) -- (-120:2*\radius) arc[start angle=-120, end angle=120, radius=2*\radius] -- cycle;

    \draw (0.35*\radius,0) arc[start angle=0, end angle=120, radius=0.35*\radius];
    \node at ({0.6*\radius*cos(60)}, {0.6*\radius*sin(60)}) {\(\theta\)};

    \draw[fill=white, draw=firstcol, very thick] (0,0) -- (120:2*\radius);
    \draw[fill=white, draw=firstcol, very thick] (0,0) -- (-120:2*\radius);

    \node [align=center] at(0.775*\legendlength, 0.95*\legendlength) {$\cup \{\infty\}$};

    \draw[-latex] (-\gridlength,0)--(\gridlength,0) node[below]{};
    \draw[-latex] (0,-\gridlength)--(0,\gridlength) node[left]{};
\end{tikzpicture}
        \label{fig:anglebounded}
    }\hspace{-0.4cm}
    \subfloat[]{%
        \begin{tikzpicture}[baseline=(current bounding box.south)]
    \pgfkeys{/pgf/number format/.cd,frac,frac whole=false}
    \pgfmathsetmacro\radius{0.95};
    \pgfmathsetmacro\angle{120};

    \pgfmathsetmacro\boxlength{1.4*\radius}
    \pgfmathsetmacro\gridlength{1.4*\radius}
    \pgfmathsetmacro\legendlength{1.25*\radius}

    \useasboundingbox (-\boxlength,-\boxlength) rectangle (\boxlength,\boxlength);
    \clip (-\boxlength,-\boxlength) rectangle (\boxlength,\boxlength);

    \begin{scope}[shift={(0.85*\radius, 0)}]
        \draw[fill=firstcol, draw=white, very thick, fill opacity = 0.2, draw opacity = 0] (0,0) -- (-120:2*\radius) arc[start angle=-120, end angle=120, radius=2*\radius] -- cycle;

        \draw (0.35*\radius,0) arc[start angle=0, end angle=120, radius=0.35*\radius];
        \node at ({0.6*\radius*cos(60)}, {0.6*\radius*sin(60)}) {\(\theta\)};

        \draw[fill=white, draw=firstcol, very thick] (0,0) -- (120:2*\radius);
        \draw[fill=white, draw=firstcol, very thick] (0,0) -- (-120:2*\radius);
    \end{scope}

    \draw[fill=black] (0.85*\radius,0) circle (0.05) node[below right=0cm and -0.05cm] {$\alpha$};

    \draw[fill=black] (-0.3*\radius,0) circle (0.05) node[above] {$\nicefrac{1}{2\rho}$};

    \draw (-0.3*\radius,0) circle (0.85*\radius);
    \begin{scope}[shift={(-0.3*\radius, 0)}]
        \draw[->, thin] (0,0) -- ++(-140:0.85*\radius) node[pos=0.5, sloped, below] {\(r\)};
    \end{scope}

    \node [align=center] at(0.775*\legendlength, 0.95*\legendlength) {$\cup \{\infty\}$};

    \draw[-latex] (-\gridlength,0)--(\gridlength,0) node[below]{};
    \draw[-latex] (0,-\gridlength)--(0,\gridlength) node[left]{};
\end{tikzpicture}
        \label{fig:semiangleproof}
    }\hspace{-0.4cm}
    \subfloat[]{%
        \begin{tikzpicture}[baseline=(current bounding box.south)]
    \pgfkeys{/pgf/number format/.cd,frac,frac whole=false}
    \pgfmathsetmacro\radius{0.95};
    \pgfmathsetmacro\angle{120};
    \pgfmathsetmacro\alphanum{0.85*\radius};

    \pgfmathsetmacro\boxlength{1.4*\radius}
    \pgfmathsetmacro\gridlength{1.4*\radius}
    \pgfmathsetmacro\legendlength{1.25*\radius}

    \useasboundingbox (-\boxlength,-\boxlength) rectangle (\boxlength,\boxlength);
    \clip (-\boxlength,-\boxlength) rectangle (\boxlength,\boxlength);

    \begin{scope}[shift={(\alphanum, 0)}]
        \draw[fill=firstcol, draw=white, very thick, fill opacity = 0.2, draw opacity = 0] (0,0) -- (-120:2*\radius) arc[start angle=-120, end angle=120, radius=2*\radius] -- cycle;

        \draw (0.35*\radius,0) arc[start angle=0, end angle=120, radius=0.35*\radius];
        \node at ({0.6*\radius*cos(60)}, {0.6*\radius*sin(60)}) {\(\theta\)};

        \draw[fill=white, draw=firstcol, very thick] (0,0) -- (120:2*\radius);
        \draw[fill=white, draw=firstcol, very thick] (0,0) -- (-120:2*\radius);
    \end{scope}

    \draw[fill=black] (\alphanum,0) circle (0.05) node[below right=0cm and -0.05cm] {$\alpha$};

    \newcommand{\drawcirclegivenrho}[2]{%
        \draw[fill=black!#2!white] ({1/(2*#1)},0) circle (0.05);
        \draw[draw=black!#2!white] ({1/(2*#1)},0) circle ({(1-2*\alphanum*#1)*0.866/(-2*#1)});
    }

    \node [align=center] at(0.775*\legendlength, 0.95*\legendlength) {$\cup \{\infty\}$};
    
    \draw[-latex] (-\gridlength,0)--(\gridlength,0) node[below]{};
    \draw[-latex] (0,-\gridlength)--(0,\gridlength) node[left]{};

    \drawcirclegivenrho{-3.33/\radius}{100}
    \drawcirclegivenrho{-1.66/\radius}{75}
    \drawcirclegivenrho{-1.1066/\radius}{50}
    \drawcirclegivenrho{-0.833/\radius}{25}
\end{tikzpicture}
        \label{fig:semiangle:non-unique}
    }\hspace{-0.2cm}
    \caption[
        The SRG of the class of $\theta$-angle-bounded operators $\operatorclass{B}_\theta$ and the construction to relate $\operatorclass{B}_\theta + \alpha \id$ to $\operatorclass{S}_{\mu,\rho}$.
    ]{
        The SRG of the class of $\theta$-angle-bounded operators $\operatorclass{B}_\theta$ and the construction to relate $\operatorclass{B}_\theta + \alpha \id$ to $\operatorclass{S}_{\mu,\rho}$.
        (a) The SRG of $\operatorclass{B}_\theta$.
        (b) The SRG of $\operatorclass{B}_\theta + \alpha \id$ and the open disk 
        $D(\nicefrac1{2\rho}, \nicefrac{\sqrt{1-4\mu\rho}}{2|\rho|})$
        from the proof of \Cref{prop:anglesem} for some $\rho < 0$ and $\mu > 0$.
        (c) The SRG of $\operatorclass{B}_\theta + \alpha \id$ and open disks
        $D(\nicefrac1{2\rho}, \nicefrac{\sqrt{1-4\mu\rho}}{2|\rho|})$ for varying $\rho < 0$, where $\mu$ is defined as in \eqref{SRG:eq:optimal-mu}.
    }
\end{figure}

\section{Semimonotonicity of circuit elements}\label{sec:circuit-elems}
	This section examines the semimonotonicity of static nonlinearities and transistors, which will then be applied to common-emitter amplifier circuits in \Cref{sec:circuits}.

\begin{figure*}
    \captionsetup[sub]{labelformat=simple}
    \centering
    \subfloat[NPN transistor]{%
        \begin{circuitikz}[american]
            \draw (0,0) node[npn, rotate=90] (p) {};
            \draw (p.B) to[short, -*] ++(0, -1.2) coordinate (B) to (B -| p.C) to[short, i^<=$i_1$] ++(-0.6, 0) coordinate (A) to[open, o-o, v^>=$v_1$] (A |- p.C) to[short, i<=$i_1$] (p.C);
            \draw (B) to (B -| p.E) to[short, i_<=$i_2$] ++(0.6,0) coordinate (C) to[open, o-o, v^>=$v_2$] (C |- p.C) to[short,i_<=$i_2$] (p.E);
        \end{circuitikz}
    } \hfill
    \subfloat[Ebers--Moll model]{%
        \begin{circuitikz}[american]
            \draw (0, 0) to[short,i_=$i_1$] ++(0.6, 0) to[short,-*] ++(2., 0) coordinate(E) to[short,*-*] ++(0, 2) coordinate(C) -- ++(-0.5, 0) coordinate(B) -- ++(0, -0.9) to[/tikz/circuitikz/bipoles/length=1.25cm, full diode, l=$i_{3}$] ++(-1.5, 0) to ++(0, 0.9) to[short,i_=$i_1$]  ++(-0.6, 0) coordinate(A);
            \draw (0, 0) to[open, o-o, v=$v_1$] (A);
            \draw (B) to ++(0, 0.0) to[/tikz/circuitikz/bipoles/length=1.25cm, cI, invert, l_=$\alpha_R i_{4}$] ++(-1.5, 0) to ++(0, -0.0);
            \draw (C) to ++(0.5, 0) coordinate(D) to ++(0, -0.9) to[/tikz/circuitikz/bipoles/length=1.25cm, full diode, l_=$i_{4}$] ++(1.5, 0) to ++(0, 0.9) to[short,i^=$i_2$] ++ (0.6, 0) coordinate (F);
            \draw (D) to ++(0, 0.0) to[/tikz/circuitikz/bipoles/length=1.25cm, cI, invert, l^=$\alpha_F i_{3}$] ++(1.5, 0) to ++(0, -0.0); 
            \draw (E) to ++(2., 0) to[short,i<_=$i_2$] ++(0.6, 0) to[open, o-o, v=$v_2$] (F);
        \end{circuitikz}
    } \hfill
    \subfloat[$\alpha_R = \frac{3}{10}, \alpha_F = \frac{2}{10}$]{%
        \begin{tikzpicture}
            \begin{axis}[
                axis line on top,
                xmin=-1, xmax=1,
                ymin=-1, ymax=1,
                axis lines=middle,
                axis line style={-stealth, thick},
                xtick=\empty,
                ytick=\empty,  
                height=2.8cm,
                scale only axis,
            ]
                \addplot[
                    only marks,
                    mark size=0.25pt,
                    color=firstcol,
                    table/col sep=comma,
                    opacity = 0.75
                ] 
                table[x=x, y=y] {data/srg_transistor2.csv};
            \end{axis}
        \end{tikzpicture}
    } \hfill
    \subfloat[$\alpha_R = \frac{110}{111}, \alpha_F = \frac{10}{11}$]{%
        \begin{tikzpicture}
            \begin{axis}[
                axis line on top,
                xmin=-1, xmax=1,
                ymin=-1, ymax=1,
                axis lines=middle,
                axis line style={-stealth, thick},
                xtick=\empty,
                ytick=\empty,  
                height=2.8cm,
                scale only axis,
            ]
                \addplot[
                    only marks,
                    mark size=0.25pt,
                    color=firstcol,
                    table/col sep=comma,
                    opacity = 0.75
                ] 
                table[x=x, y=y] {data/srg_transistor.csv};
            \end{axis}
        \end{tikzpicture}
    }
    \caption{NPN transistor.
    (a) Two-port model.
    (b) Ebers--Moll model.
    (c)-(d) Numerical SRG of the NPN transistor $G(T_{\rm NPN})$ for different values of $\alpha_R$ and $\alpha_F$.
    In both cases, the incremental angle is upper bounded by $135$ degrees.
    }
    \vspace{-12pt}
    \label{fig:transebers}
\end{figure*}

\subsection{One-dimensional static nonlinearities}
\label{subsec:tunnel}
Static nonlinearities are simple mathematical models which can be used to describe behaviors of certain circuit elements, such as nonlinear resistances and tunnel diodes.
The following result, which can also be derived from \cite[Prop.\ 9]{chaffey2023graphical} and \cite[Prop.\ 4.14(iv)]{evens2023convergence}, shows that any one-dimensional static nonlinearity with a bounded slope is semimonotone. For clarity, we provide a self-contained proof here.

\begin{proposition}\label{cor:semi:SISO}
    Let $T : \Re \rightarrow \Re$ be a one-dimensional single-valued operator and let $\ell > 0$ and $\sigma \in (-\ell, \ell]$.
    Then, the slope of $T$ is bounded between $\sigma$ and $\ell$, i.e., 
    \begin{align}
        \label{eq:semi:SISO:cond}
        \sigma \leq \frac{T(x) - T(y)}{x - y} \leq \ell, 
        \qquad  \forall x, y \in \Re, x \neq y,
    \end{align}
    if and only if $T$ is
    \(
        \bigl(
            \frac{\sigma\ell}{\ell + \sigma},
            \frac{1}{\ell + \sigma}
        \bigr)
    \)-semimonotone.
\end{proposition}
\begin{proof}
    Since $T : \Re \rightarrow \Re$ is
    one-dimensional and single-valued, it holds by \Cref{def:SRG} that 
    \begin{align*}
        G(T)
            =
        \left\{
            \frac{T(x)-T(y)}{x-y} \mathrel{\bigg|\;x, y \in \Re, x \neq y}
        \right\}.
    \end{align*}
    As the class of semimonotone operators is SRG-full, it follows from \Cref{prop:srgsemi,prop:semisrgfull} that $T$ is
    \(
        \bigl(
            \frac{\sigma\ell}{\ell + \sigma},
            \frac{1}{\ell + \sigma}
        \bigr)
    \)-semimonotone
    if and only if
    \(
        G(T)
        \subseteq
        \closure D(\tfrac{\ell + \sigma}{2}, \tfrac{\ell - \sigma}{2}).
    \)
    Since this is equivalent to \eqref{eq:semi:SISO:cond}, the proof is completed.
\end{proof}
Note that \eqref{eq:semi:SISO:cond} is closely related to the notion of sector nonlinearities, which satisfy
\(
    \sigma \leq \frac{T(x)}{x} \leq \ell
\)
for all $x \in \Re$.
To highlight the importance of this result, consider its application to tunnel diodes: despite having a \emph{negative resistance region}, these elements remain semimonotone as long as their slope is bounded, as illustrated below.
Tunnel diodes have long been used to model oscillatory behavior in dynamical systems, including the van der Pol oscillator and the FitzHugh--Nagumo neuron model \cite{nagumo1962active}.
\begin{example}
    \label{exmp:tunnel_diode}
    Let $r_1,r_2,\bar{v} > 0$ with $r_1 < r_2$ and consider the piecewise linear tunnel diode 
    \[
        T_{\rm tunnel} : \mathbb{R} \to \mathbb{R} : v \mapsto \begin{cases}
            r_1^{-1}(v + \bar{v}) +r_2^{-1}\bar{v}, \quad& \textnormal{if $v < -\bar{v}$,}\\
            -r_2^{-1}v, \quad& \textnormal{if $|v| \leq \bar{v}$,} \\
            r_1^{-1}(v-\bar{v}) - r_2^{-1}\bar{v}, \quad& \textnormal{if $v > \bar{v}$.}
        \end{cases}
    \]

    By applying \Cref{cor:semi:SISO} with $\sigma = -r_2^{-1}$ and $\ell = r_1^{-1}$, it follows that $T_{\rm tunnel}$ is 
    \(
    \bigl(
        \frac{1}{r_1-r_2}, \frac{r_1r_2}{r_2-r_1}
    \bigr)
    \)-semimonotone.
\end{example}

We note that a non-tight variant of \Cref{cor:semi:SISO} was used in \cite[Ex.~3]{huijzer2025modelling} to show that memristors admitting both flux- and charge-controlled representations are semimonotone with positive moduli.

\subsection{Transistor}\label{subsec:transistor}
Let 
$
    \matr R \coloneqq
    \begin{psmallmatrix}
        1 & -\alpha_R  \\
        - \alpha_F & 1
    \end{psmallmatrix}
$.
An NPN transistor \cite[Ch.\ 17.3]{desoer1969basic} can then be represented using the Ebers--Moll model
\begin{equation}
    T_{\rm NPN}
    \begin{pmatrix}
        v_{1} \\ v_{2}
    \end{pmatrix}
        {}\coloneqq{}
    \left\{
        R
        \begin{pmatrix}
            u_1 \\ u_2
        \end{pmatrix}
        \mathrel{\bigg|}
        \begin{aligned}
            &u_1 {}\in{} T_{\rm D}(v_1)\\
            &u_2 {}\in{} T_{\rm D}(v_2) 
        \end{aligned}
    \right\},
    \label{eq:ebers-moll}
\end{equation}
where $T_{\rm D}$ is the device law of a diode (which is typically monotone), $\alpha_R$ is the reverse short-circuit current ratio and $\alpha_F$ is the forward short-circuit current ratio. These current ratios are usually around $0.9$ to $0.995$ \cite[p.\ 725]{desoer1969basic}. \Cref{fig:transebers} shows the equivalent circuit representation of this model, as well as a numerical computation of its SRG by sampling random points, which suggests that this operator is angle-bounded. We provide a formal proof of this observation in \Cref{prop:transistor:angle-bounded}, under the following classical assumption.

\begin{assumption}
    \label{ass:NPN}
    $T_{\rm NPN}$ is given by \eqref{eq:ebers-moll}, with $\alpha_R, \alpha_F \in [0,1)$ and $T_{\rm D} \in \operatorclass{M}$ is monotone and outer semicontinuous.
\end{assumption}

Note that outer semicontinuity of $T_{\rm D}$ ensures by definition that also $T_{\rm NPN}$ is also outer semicontinuous.

\begin{proposition}\label{prop:transistor:angle-bounded}
    Suppose that \Cref{ass:NPN} holds.
    Then, $T_{\rm NPN}$ is $\theta$-angle-bounded, where
    \[
        \theta
            {}\coloneqq{}
        \frac{\pi}{2} + \max(\arctan(\alpha_F), \arctan(\alpha_R)).
    \]
\end{proposition}
\begin{proof}
    $
        T_{\rm NPN} = \matr R \circ (T_{\rm D} \times T_{\rm D})
    $
    is $\theta$-angle-bounded if and only if
    \begin{align*}
        \angle \left(x - y, \matr R (u - v)\right) \leq \theta,
        \;\;
        \forall (x, u), (y, v) \in \graph (T_{\rm D} \times T_{\rm D}).
    \end{align*}
    Defining $z \coloneqq x-y$ and $w \coloneqq u-v$, this condition holds for any monotone operator $T_{\rm D}$ if and only if
    \begin{align}
        \angle \left(z, \matr R (w)\right) \leq \theta,
        \;\;
        \forall z, w \in \Re^2 : z_1 w_1 \geq 0, z_2 w_2 \geq 0.
        \label{eq:proof:thm:transistor:cond}
    \end{align}
    If $\|z\| = 0$ or $\|w\|=0$ then \eqref{eq:proof:thm:transistor:cond} holds vacuously.
    Otherwise, let $\phi \coloneqq \atantwo(z_2, z_1)$ and $\psi \coloneqq \atantwo(w_2, w_1)$, so that
    \begin{equation*}
        z = \|z\|\begin{bmatrix}
            \cos(\phi)\\
            \sin(\phi)
        \end{bmatrix}
        \quad\text{and}\quad
        w = \|w\|\begin{bmatrix}
            \cos(\psi)\\
            \sin(\psi)
        \end{bmatrix}.
    \end{equation*}
    Defining $\delta_k \coloneqq [-\pi + \frac{\pi}{2}k,-\frac{\pi}{2} + \frac{\pi}{2}k]$, it holds that
    \(
        z_1 w_1 \geq 0
    \)
    and 
    \(
        z_2 w_2 \geq 0
    \)
    if and only if there exists a $k \in \{0, 1, 2, 3\}$ such that 
    \(
        \phi, \psi \in \delta_k.
    \)
    Therefore, by definition of $\angle$, \eqref{eq:proof:thm:transistor:cond} holds if and only if
    \begin{align}
        \arccos\left(\frac{f(\phi, \psi)}{g(\phi, \psi)}\right)
            \leq
        \theta,
        \;\;
        \forall \phi, \psi \in \delta_k, k \in \{0, 1, 2, 3\},
        \label{eq:proof:thm:transistor:arccoscond}
    \end{align}
    where
    \begin{align*}
        f(\phi, \psi)
            \coloneqq{}&
        \left\langle 
        \begin{bmatrix}
            \cos(\phi)\\
            \sin(\phi)
        \end{bmatrix}, R\begin{bmatrix}
            \cos(\psi)\\
            \sin(\psi)
        \end{bmatrix}\right\rangle,\\
        g(\phi, \psi)
            \coloneqq{}&
        \left\|
            R\begin{bmatrix}
                \cos(\psi)\\
                \sin(\psi)
            \end{bmatrix}
        \right\| 
            >
        0,
    \end{align*}
    so that
    \begin{equation*}
        \resizebox{\linewidth}{!} 
{       $
        \frac{f(\phi, \psi)}{g(\phi, \psi)}
            ={}
        \frac{
            \cos(\phi - \psi) - \alpha_R \cos(\phi)\sin(\psi) - \alpha_F \sin(\phi)\cos(\psi)
        }
        {
            \sqrt{1 + \alpha_R^2 \sin^2(\psi) - 2(\alpha_R + \alpha_F)\sin(\psi)\cos(\psi) + \alpha_F^2\cos^2(\psi)}
        }.$}
    \end{equation*}
    If $\phi, \psi \in \delta_1$ or $\phi, \psi \in \delta_3$, then $f(\phi, \psi)$ is lower bounded by zero, and consequently so is $\tfrac{f(\phi, \psi)}{g(\phi, \psi)}$.
    If $\phi, \psi \in \delta_0$ or $\phi, \psi \in \delta_2$, then 
    \[
        \frac{f(\phi, \psi)}{g(\phi, \psi)}
            \geq
        \min\Biggl\{
            \frac{-\alpha_F}{\sqrt{1+\alpha_F^2}},
            \frac{-\alpha_R}{\sqrt{1+\alpha_R^2}}
        \Biggr\}.
    \]
    Note that this bound is tight and attained for $(\phi, \psi) = (\nicefrac{\pi}{2}, 0)$ if $\alpha_F \geq \alpha_R$ and for $(0, \nicefrac{\pi}{2})$ otherwise.
    Therefore, the claim follows directly from \eqref{eq:proof:thm:transistor:arccoscond}
    by observing that $\arccos$ is a decreasing function and that for any $\alpha \in \Re$
    \[
        \arccos
        \left(
            \frac{-\alpha}{\sqrt{1+\alpha^2}}
        \right)
            =
        \frac{\pi}{2} + \arctan(\alpha).\qedhere
    \]
\end{proof}

\begin{cor}\label{cor:transistor:semi}
    Suppose that \Cref{ass:NPN} holds.
    Then, $T_{\rm NPN}$ is $(\frac{1}{8\rho}, \rho)$-semimonotone for all $\rho < 0$.
\end{cor}
\begin{proof}
    By \Cref{prop:transistor:angle-bounded}, $T_{\rm NPN}$ is $\frac{3\pi}{4}$-angle-bounded for any $\alpha_F, \alpha_R \in [0,1)$. The claim then follows immediately from \Cref{prop:anglesem}.
\end{proof}

\section{Nonmonotone common-emitter amplifiers}\label{sec:circuits}
	In this section, we consider common-emitter amplifiers involving transistors and tunnel diodes, and show that their response can be computed efficiently using Chambolle--Pock.

\subsection{Background on Chambolle--Pock for circuit theory}
For many circuits, standard methods of loop and cut-set analysis can be used to find a hybrid representation of the circuit~\cite{desoer1969basic}.
For instance, denote the internal currents and voltages by respectively $i \in \mathcal{H}^n$ and $v \in \mathcal{H}^m$ and let the influence of the current and voltage sources on the network be given by respectively $s_i \in \mathcal{H}^{n}$ and $s_v \in \mathcal{H}^{m}$. Then, the behavior of a circuit can be retrieved by solving an inclusion problem of the form
\begin{align}\label{eq:circuits:inclusion}
    0\in
    \begin{bmatrix}
        R(i)\\ G(v)
    \end{bmatrix}
        +
    \begin{bmatrix}
        0 & L^\top\\ 
        - L & 0
    \end{bmatrix}  
    \begin{bmatrix}
        i\\ 
        v
    \end{bmatrix}
        +
    \begin{bmatrix}
        s_v\\ 
        s_i
    \end{bmatrix}
\end{align}
where $R : \mathcal{H}^n \rightrightarrows \mathcal{H}^n$ and $G : \mathcal{H}^m \rightrightarrows \mathcal{H}^m$ are operators containing respectively the resistive and conductive elements present in the network and the matrix $L \in \Re^{m \times n}$ encodes Kirchhoff's current and voltage laws.
When constructing a hybrid representation, series or parallel interconnections of semimonotone operators can be grouped using basic calculus rules for the sum and parallel sum of semimonotone operators; see \cite[Prop.\ 4.7]{evens2023convergence}. Note that these rules can also be derived graphically using the tools introduced in \Cref{sec:srg-nonmon}.

A method that is particularly suited for solving inclusion problems of this form is the Chambolle--Pock algorithm (CPA)~\cite{chambolle2011firstorder} (also known as the primal-dual hybrid gradient (PDHG) method \cite{zhu2008efficient}).
When applied to \eqref{eq:circuits:inclusion}, this algorithm performs alternating updates of internal currents and voltages using resolvent computations, where the $\gamma$-resolvent of an operator $T : \mathcal{H} \rightrightarrows \mathcal{H}$ is defined as $J_{\gamma T} \coloneqq (\id + \gamma T)^{-1}$.
In particular, for strictly positive stepsizes \(\gamma,\tau > 0\), relaxation parameter $\lambda > 0$ and an initial guess \((i^0,v^0)\in\mathcal{H}^n \times \mathcal{H}^m\), the CPA is given by
\begin{equation}\tag{CPA}\label{eq:CP}
    \begin{aligned}
        \bar i^{k}
            {}\in{}&
        J_{\gamma \tilde R}\big(i^k - \gamma L^\top v^k\big)\\
        \bar v^{k}
            {}\in{}&
        J_{\tau \tilde G}\big(v^k + \tau L(2\bar{i}^{k} - i^k)\big)\\
        i^{k+1}
            {}\in{}&
        i^k + \lambda (\bar i^k - i^k)\\
        v^{k+1}
            {}\in{}&
        v^k + \lambda (\bar v^k - v^k)
    \end{aligned}
\end{equation}
where $\tilde R(i) \coloneqq R(i) + s_v$ and $\tilde G(v) \coloneqq G(v) + s_i$.
The idea of applying \eqref{eq:CP} for solving \eqref{eq:circuits:inclusion} was previously explored in~\cite{chaffey2023circuit} in the monotone setting.
One of the main advantages of this methodology, besides its scalability, is that it allows $R$ and $G$ to be multi-valued, which standard Newton solvers for electrical circuits are unable to deal with. 

This same challenge also motivates the nonsmooth dynamics framework developed in~\cite{acary_nonsmooth_2011}, which addresses multi-valuedness in analog switched circuits through the use of complementarity problems and inclusions involving normal cones.

\begin{figure}
    \centering
    \subfloat[Leakage current model]{%
        \begin{circuitikz}[american]
            \draw (0,0) node[npn, rotate=90] (p) {};
            \draw (p.B) to[short, -*] (0, -2) coordinate (B) to (B -| p.C) to[short, i^<=$i_1$] ++(-0.6, 0) coordinate (A) to[open, o-o, v^>=$v_1$] (A |- p.C) to[short, i<=$i_1$] (p.C);
            \draw (B) to (B -| p.E) to[short, i_<=$i_2$] ++(0.6,0) coordinate (C) to[open, o-o, v^>=$v_2$] (C |- p.C) to[short,i_<=$i_2$] (p.E);
            \draw (p.E) to[/tikz/circuitikz/bipoles/length=1.25cm, R, l_=$r$] (p.E |- C);
            \draw (p.C) to[/tikz/circuitikz/bipoles/length=1.25cm, R, l^=$r$] (p.C |- A);
        \end{circuitikz}
        \label{fig:NPN-leaky}
    } \hfill
    \subfloat[Amplifier model]{%
        \begin{circuitikz}[american]
            \ctikzset{bipoles/resistor/height=0.1}
            \ctikzset{bipoles/resistor/width=0.3}
            \ctikzset{bipoles/ageneric/height=0.2}
            \ctikzset{bipoles/ageneric/width=0.45}
    
            \draw (0,0) node[npn] (p) {};
            \draw (p.E) to[R, -, l^=$R_E$] ++(0, -0.65) to[short, -] ++(0, -0.2) coordinate (G) to[short] ++(-1.35, 0) coordinate (A) to[/tikz/circuitikz/bipoles/length=1cm, V, invert, l=$v_{\textnormal{in}}$] (A |- p.B) to (p.B);
    
            \draw (p.B) to ++(0.05, 0) to ++(0, -0.47) to[R, bipoles/resistor/height=0.1, bipoles/resistor/width=0.3, l_=$r$] ++(0.8, 0);
            \draw (p.B) to ++(0.05, 0) to ++(0, 0.47) to[R, bipoles/resistor/height=0.1, bipoles/resistor/width=0.3, l^=$r$] ++(0.8, 0);
    
            \draw (p.C) to[R, l_=$R_C$] ++(1.5, 0) coordinate (B) to[/tikz/circuitikz/bipoles/length=1cm, V, l= $v_+$] (B |- A) to ++(A |- p.B) node[ground]{} to (G);
            \draw (0.5, 0.65) rectangle ++(0.6, 0.25);
            \draw[fill] (1.02, 0.65) rectangle ++(0.08, 0.25);

            \draw (-0.13, -0.83) rectangle ++(0.25, -0.6);
            \draw[fill] (-0.13, -1.35) rectangle ++(0.25, -0.08);

            \node at (0.3, -0.7) {$+$};
            \node at (0.3, -1.45) {$-$};

            \node at (0.3, 0.45) {$+$};
            \node at (1.2, 0.45) {$-$};            
        \end{circuitikz}
        \label{fig:common-emitter amplifier}
    }
    \caption{Common-emitter amplifier with leakage current.}
\end{figure}

\subsection{Implementation details}
In all upcoming numerical examples, the transistor parameters are given by $\alpha_R = \frac{110}{111}, \alpha_F = \frac{10}{11}$ and the internal diodes are modeled by ideal diodes, i.e., by
\[
    T_{\rm D}:
    v
        \mapsto
    \begin{cases}
        \{0\}, \quad& \text{if } v < 0,\\
        [0, +\infty), \quad& \text{if } v = 0.  
    \end{cases}
\]
We denote the transistor voltages by $v_1$ and $v_2$.
For numerical simulations, the iterations are stopped once the norm of the relative difference between successive iterates is less than $\epsilon = 10^{-8}$.
We do not explicitly verify that the resolvent has full domain, as this assumption mainly ensures the global well-definedness of iterates.
No such domain-related issues were encountered in any of our simulations. The code for reproducing the experiments is publicly available\footnote{\url{https://github.com/JanQ/SRG-networks}}.

\subsection{Examples}
First, consider the NPN transistor shown in \Cref{fig:NPN-leaky}, which includes two leakage resistors with resistance $r$.
For this model, the voltages $v = (v_{1}, v_{2})$ which yield a desired current $i = (i_1, i_2)$ can be obtained by solving the following inclusion problem:
\begin{equation}
    \begin{pmatrix}
        i_1\\
        i_2
    \end{pmatrix}
        {}\in{}
    T_{\rm NPN}
    \begin{pmatrix}
        v_{1}\\
        v_{2}
    \end{pmatrix}
        +
    \frac{1}{r}
    \begin{pmatrix}
        v_{1}\\
        v_{2}
    \end{pmatrix}.
    \label{eq:inclusion:leaky-NPN}
\end{equation}

As established in \Cref{prop:transistor:angle-bounded}, an NPN transistor is $\frac{3\pi}{4}$-angle-bounded under standard assumptions summarized in \Cref{ass:NPN}.
Therefore, for any $r>0$ it holds that $T_{\rm NPN, r} \coloneqq T_{\rm NPN} + \frac{1}{r}\id$ is 
\begin{enumerate}[label=(\roman*)]
    \item $\tfrac{r(1 - \sqrt2)}{2}$-comonotone owing to \Cref{cor:anglesem:identity},
    \item $(\tfrac1{2r}, -\tfrac{r}{2})$-semimonotone owing to \Cref{prop:anglesem}.
\end{enumerate}

The following proposition demonstrates how this result can be leveraged to solve \eqref{eq:inclusion:leaky-NPN} directly using the proximal point algorithm (PPA), using known convergence results from \cite[Thm.\ 2.4, Tab.\ 1]{evens2023convergence} in the comonotone setting.
This is verified numerically in \Cref{fig:leakytransistorppa}.
\begin{proposition}
    Suppose that \Cref{ass:NPN} holds and let $r > 0$.
    Define $\tilde T_{\rm NPN, r}(v) \dfn T_{\rm NPN, r}(v) -  i$.
    Then, any sequence $(v^k)_{k \in \N}$ satisfying the PPA update rule
    \[
        v^{k+1}
            \in 
        J_{\gamma \tilde T_{\rm NPN, r}}(v^k)
    \]
    with stepsize $\gamma > r(\sqrt{2} - 1)$ converges to a solution of \eqref{eq:inclusion:leaky-NPN}.
\end{proposition}
\begin{proof}
    Follows directly from \cite[Thm.\ 2.4]{evens2023convergence}, since $\tilde T_{\rm NPN, r}$ is outer semicontinuous and $\tfrac{r(1 - \sqrt2)}{2}$-comonotone.
\end{proof}

\begin{figure}
    \centering
    \captionsetup[subfigure]{labelformat=empty}
    \subfloat[]{%
        \begin{tikzpicture}
            \begin{axis}[
                xlabel={$t$ [\si{\second}]},
                xlabel style={yshift=0ex, font=\footnotesize},
                xticklabel style={font=\footnotesize},
                ylabel={$i$ [\si{\ampere}]},
                ylabel style={yshift=-1ex, font=\footnotesize},
                yticklabel style={font=\footnotesize},
                legend pos=south west,
                legend image code/.code={\draw[#1] (0,0) -- (3.2mm,0);},
                legend style={font=\footnotesize},
                height=2.65cm,
                xmin=0, xmax=1,
                ymin=-1.175, ymax=1.175,
                scale only axis,
                clip mode=individual
            ]
            \addplot[color=firstcol, thick] 
                table[x=t, y=i1, col sep=comma] {data/transistor_leaky_ppa.csv};
            \addlegendentry{$i_1$} 

            \addplot[color=secondcol,densely dashed, thick] 
                table[x=t,y=i2, col sep=comma] {data/transistor_leaky_ppa.csv};
            \addlegendentry{$i_2$} 
            \end{axis}
        \end{tikzpicture}
    } \hfill
    \subfloat[]{%
        \begin{tikzpicture}
            \begin{axis}[
                xlabel={$t$ [\si{\second}]},
                xlabel style={yshift=0ex, font=\footnotesize},
                xticklabel style={font=\footnotesize},
                ylabel={$v$ [\si{\volt}]},
                ylabel style={yshift=-1ex, font=\footnotesize},
                yticklabel style={font=\footnotesize},
                legend pos=south west,
                legend image code/.code={\draw[#1] (0,0) -- (3.2mm,0);},
                legend style={font=\footnotesize},
                height=2.65cm,
                scale only axis,
                xmin=0, xmax=1,
                ymin=-10, ymax=2,
                clip mode=individual
            ]
            \addplot[color=firstcol, thick] 
                table[x=t, y=v1, col sep=comma] {data/transistor_leaky_ppa.csv};
            \addlegendentry{$v_1$} 

            \addplot[color=secondcol,densely dashed, thick] 
                table[x=t,y=v2, col sep=comma] {data/transistor_leaky_ppa.csv};
            \addlegendentry{$v_2$} 
            \end{axis}
        \end{tikzpicture}
    }
    \vspace{-20pt}
    \caption{
        Solution of inclusion problem \eqref{eq:inclusion:leaky-NPN} with leakage resistance $r=\SI{10}{\ohm}$ for a given desired sinusoidal current $i$. The solution has been obtained after $27$ proximal point iterations with constant stepsize $\gamma = r > r(\sqrt{2} - 1)$.
    } \label{fig:leakytransistorppa}
\end{figure}

In what follows, we connect the nonmonotone element $T_{\rm NPN, r}$ in a so-called common-emitter amplifier circuit, visualized by \Cref{fig:common-emitter amplifier}.
Denote the voltage sources by $v_{\rm in}\in\mathbb{R}$ and $v_+ \in \mathbb{R}_{++}$. With a suitable convention of the directions of the internal currents and voltages, the response of this network is given by \eqref{eq:circuits:inclusion}, where 
\begin{equation}\label{eq:common-emitter:vars}
    \begin{aligned}
        R
            \coloneqq{}&
        R_C \times R_E,
        \quad
        G
            \coloneqq
        T_{\rm NPN,r},
        \quad
        L
            \coloneqq
        \matr{I}_2,\\
        \quad
        s_v
            \coloneqq{}&
        \begin{pmatrix}
            v_+ - v_{\textnormal{in}} \\ -v_{\textnormal{in}}
        \end{pmatrix},
        \quad
        s_i
            \coloneqq
        0.
    \end{aligned}
\end{equation}

By leveraging the convergence results for \eqref{eq:CP} from \cite{evens2023convergenceCP}, we show in the following proposition that the response of this common-emitter amplifier circuit can be computed by \eqref{eq:CP} under suitable conditions on $R_E$ and $R_C$ and corresponding stepsize conditions.
Specifically, \Cref{ex:common-emitter:str-mon} considers the setting where the resistances are strongly monotone, for which a numerical example is provided in \Cref{fig:commonemittercpa}, while \Cref{ex:common-emitter:semi} considers the setting where they are merely $(\tfrac{9}{8}r, -\tfrac1{8r})$-semimonotone. This latter setting includes resistances which have a \emph{negative resistance region} and is illustrated numerically in \Cref{fig:commonemittertunnelcpa}.

\begin{figure}[t]
    \centering
    \captionsetup[subfigure]{labelformat=empty}
    \subfloat[]{%
        \begin{tikzpicture}
            \begin{axis}[
                xlabel style={yshift=0ex, font=\footnotesize},
                xticklabel style={font=\footnotesize},
                ylabel style={yshift=-3ex, font=\footnotesize},
                yticklabel style={font=\footnotesize},
                legend image code/.code={\draw[#1] (0,0) -- (3.2mm,0);},
                legend style={font=\footnotesize},
                height=2.65cm,
                scale only axis,
                tick scale binop=\times,
                xlabel={$t$ [\si{\second}]},
                ylabel={$i$ [\si{\ampere}]},
                xmin=0, xmax=2,
                ymin=-0.07, ymax=0.04,
                legend pos=south west,
                yticklabel style={
                    /pgf/number format/fixed,
                    /pgf/number format/precision=3
                },
                scaled ticks = false,
                scaled y ticks = false,
            ]
            \addplot[color=firstcol, thick] 
            table[x=t, y=ic, col sep=comma] {data/common_emitter.csv};
            \addlegendentry{$i_C$} 
            \addplot[color=secondcol, densely dashed, thick] 
            table[x=t,y=ie, col sep=comma] {data/common_emitter.csv};
            \addlegendentry{$i_E$} 
            \end{axis}
        \end{tikzpicture}
    } \hfill
    \subfloat[]{%
        \begin{tikzpicture}
            \begin{axis}[
                xlabel style={yshift=0ex, font=\footnotesize},
                xticklabel style={font=\footnotesize},
                ylabel style={yshift=-1ex, font=\footnotesize},
                yticklabel style={font=\footnotesize},
                legend image code/.code={\draw[#1] (0,0) -- (3.2mm,0);},
                legend style={font=\footnotesize},
                height=2.65cm,
                scale only axis,
                tick scale binop=\times,
                legend pos=south west,
                xlabel={$t$ [\si{\second}]},
                ylabel={$v$ [\si{\volt}]},
                xmin=0, xmax=2,
                ymin=-2.5, ymax=0.2
            ]
            \addplot[color=firstcol, thick] 
                table[x=t,y=v1,col sep=comma] {data/common_emitter.csv};
            \addlegendentry{$v_1$} 
            \addplot[color=secondcol,densely dashed, thick] 
            table[x=t,y=v2, col sep=comma] {data/common_emitter.csv};
            \addlegendentry{$v_2$} 
            \end{axis}
        \end{tikzpicture}
    }
    \vspace{-20pt}
    \caption{
        Calculated internal variables $i_C,i_E, v_1, v_2$ for the common-emitter amplifier with linear resistors $R_E = r_E\id$ and $R_C = r_C\id$, where $r_E$ and $r_C$ are strictly larger than $\nicefrac{r(\sqrt{2}-1)}{2}$.
        The circuit parameters are $v_{\rm in} = \sin(2\pi t) \si{\volt}$, $v_+ =\SI{5}{\volt}$, $r_E=\SI{30}{\ohm}$, $r_C=\SI{150}{\ohm}$ and $r=\SI{100}{\ohm}$.
        These internal variables were obtained using Chambolle--Pock after $617$ iterations with stepsizes $\gamma = 0.001$, $\tau=700$ and relaxation parameter $\lambda=1$.
    } \label{fig:commonemittercpa}
\end{figure}

\begin{figure}[t]
    \centering
    \captionsetup[subfigure]{labelformat=empty}
    \subfloat[]{%
        \begin{tikzpicture}
            \begin{axis}[
                xlabel style={yshift=0ex, font=\footnotesize},
                xticklabel style={font=\footnotesize},
                ylabel style={yshift=-3ex, font=\footnotesize},
                yticklabel style={font=\footnotesize},
                legend image code/.code={\draw[#1] (0,0) -- (3.2mm,0);},
                legend style={font=\footnotesize},
                height=2.65cm,
                scale only axis,
                tick scale binop=\times,
                xlabel={$t$ [\si{\second}]},
                ylabel={$i$ [\si{\ampere}]},
                xmin=0, xmax=2,
                ymin=-0.01, ymax=0.015,
                legend pos=south west,
                samples=500,
                yticklabel style={
                    /pgf/number format/fixed,
                    /pgf/number format/precision=3
                },
                scaled ticks = false,
                scaled y ticks = false,
            ]
            \addplot[color=firstcol, thick] 
            table[x=t, y=ic, col sep=comma] {data/common_emitter_tunnel.csv};
            \addlegendentry{$i_C$} 

            \addplot[color=secondcol,densely dashed, thick] 
            table[x=t,y=ie, col sep=comma] {data/common_emitter_tunnel.csv};
            \addlegendentry{$i_E$} 
            \end{axis}
        \end{tikzpicture}
        \label{fig:subfig1}%
    } \hfill
    \subfloat[]{%
        \begin{tikzpicture}
            \begin{axis}[
                xlabel style={yshift=0ex, font=\footnotesize},
                xticklabel style={font=\footnotesize},
                ylabel style={yshift=-1ex, font=\footnotesize},
                yticklabel style={font=\footnotesize},
                ytick={-6,-5,-4},
                legend image code/.code={\draw[#1] (0,0) -- (3.2mm,0);},
                legend style={font=\footnotesize},
                height=2.65cm,
                scale only axis,
                tick scale binop=\times,
                xlabel={$t$ [\si{\second}]},
                ylabel={$v_{\rm tunnel}$ [\si{\volt}]},
                xmin=0, xmax=2,
                ymin=-6, ymax=-3.8,
                legend pos=south west
            ]
            \addplot[color=firstcol, thick] 
            table[x=t, y=vtunnel, col sep=comma] {data/common_emitter_tunnel.csv};
            \end{axis}
        \end{tikzpicture}
        \label{fig:subfig2}%
    }
    \vspace{-20pt}
    \caption{
        Calculated internal variables $i_C,i_E, v_{\rm tunnel}$ for the common-emitter amplifier with a linear resistor $R_E = r_E\id$ and a (multi-valued) resistor $R_C = T_{\rm tunnel}^{-1}$ defined as the inverse of the tunnel diode in \Cref{exmp:tunnel_diode}. The circuit parameters are $v_{\textnormal{in}} = \sin(2\pi t)\si{\volt}, v_+=\SI{5}{\volt}$, $r = \SI{100}{\ohm}$ and $r_E = \SI{100}{\ohm}$. The tunnel diode parameters are $r_1=\SI{100}{\ohm}$, $r_2=\SI{900}{\ohm}$ and $\bar{v} = \SI{5}{\volt}$. By \Cref{cor:semi:SISO}, $R_E$ and $R_C$ are $(\tfrac{900}{8}, -\tfrac1{800})$-semimonotone and the conditions of \Cref{ex:common-emitter:semi} are satisfied.
        These internal variables were obtained using Chambolle--Pock after $223$ iterations with stepsizes $\gamma = \tfrac{5}{9r} = \nicefrac1{180}$, $\tau = \tfrac{8r}{5} = 160$ and relaxation parameter $\lambda = \nicefrac14$.
        Note that both the positive and negative resistance regions of the tunnel diodes are encountered during this experiment.
    }
    \label{fig:commonemittertunnelcpa}
\end{figure}

\begin{proposition}\label{ex:common-emitter}
    Consider problem \eqref{eq:circuits:inclusion}, where $R, G, L, s_v$ and $s_i$ are defined as in 
    \eqref{eq:common-emitter:vars}.
    Suppose that \Cref{ass:NPN} holds, that the leakage resistance 
    $r > 0$ and that the (nonlinear) operators $R_C : \Re \rightarrow \Re$ and $R_E : \Re \rightarrow \Re$ in \eqref{eq:common-emitter:vars} are outer semicontinuous.
    Suppose that $R_C$ and $R_E$ are 
    \begin{enumerate}
        \item \label{ex:common-emitter:str-mon}(either) $\sigma$-monotone for some
        \(
            \sigma > \tfrac{r(\sqrt2 - 1)}{2}
        \)
        and 
        \[
            \gamma \in (0, \nicefrac{1}{\underline{\tau}}),\quad \tau \in (\underline{\tau}, \nicefrac{1}{\gamma}),\quad \lambda\in \bigl(0, 
            2(1-\nicefrac{\underline{\tau}}{\tau})
            \bigr),
        \]
        where
        \(
            \underline{\tau}
                \coloneqq
            -\tfrac{\sigma r(1 - \sqrt2)}{r(1 - \sqrt2) + 2\sigma},
        \)
        \item \label{ex:common-emitter:semi}(or) $(\tfrac{9}{8}r, -\tfrac1{8r})$-semimonotone and 
        \begin{align*}
            \gamma 
                {}\in{}&
            \Bigl(
                \tfrac{5 - \sqrt{10}}{9r},
                \tfrac{5 + \sqrt{10}}{9r}
            \Bigr), \quad
            \tau
                \in 
            \Bigl(
                \tfrac{9r(6r\gamma-1)}{51r\gamma-10},
                \tfrac{1}{\gamma}
            \Bigr),\\
            \lambda
                {}\in{}&
            \left(
                0, 
                2\left(1+\Delta_{\gamma,\tau} - \sqrt{\Delta_{\gamma,\tau}^2 - \tfrac3{20\gamma\tau}(1-\gamma\tau)}\right)
            \right),
        \end{align*}
        where
        \(
            \Delta_{\gamma,\tau}
                \coloneqq
            -\tfrac1{12r\gamma} - \tfrac{9r}{20\tau}
        \).
    \end{enumerate}
    Then, any sequence $(i^k, v^k)_{k \in \N}$ generated by \eqref{eq:CP} with stepsizes $\gamma$ and $\tau$ and relaxation parameter $\lambda$ converges to a solution of \eqref{eq:circuits:inclusion}.
\end{proposition}
\begin{proof}
    Note that both $\tilde R$ and $\tilde G^{-1}$ are outer semicontinuous.
    \begin{enumerate}[label=(\roman*)]
        \item Follows directly from \cite[Cor.\ 5.2]{evens2023convergenceCP}
        with $\beta_{\rm P} = 0$ and $\beta_{\rm D} = -\underline{\tau} < 0$,
        using that 
        $\tilde R$ is $\sigma$-monotone
        and
        $\tilde G^{-1}$ is $\tfrac{r(1 - \sqrt2)}{2}$-monotone (since $\tilde{G} = G$ is $\tfrac{r(1-\sqrt{2})}{2}$-comonotone).
        \item By assumption,
        $\tilde R$ is $(\tfrac{9}{8}r, -\tfrac1{8r})$-semimonotone
        and 
        $\tilde G^{-1}$ is $(-\tfrac{r}{2}, \tfrac1{2r})$-semimonotone.
        The claim then follows from \cite[Cor.\ 5.2]{evens2023convergenceCP} with $\beta_{\rm P} = -\tfrac{1}{6r}$ and $\beta_{\rm D} = -\tfrac{9}{10}r$, where we used that 
        $\theta_{\gamma\tau}(\|L\|) = \sqrt{\Delta_{\gamma,\tau}^2 - \tfrac3{20\gamma\tau}(1-\gamma\tau)}$ as defined in \cite[Eqns.\ (12) \& (20)]{evens2023convergenceCP} (the definition of $\theta_{\gamma\tau}$ in \cite[Relaxation parameter rule II]{evens2023convergenceCP} contains a typographical error).
        \qedhere
    \end{enumerate}
\end{proof}

\section{Conclusion}
	In this work, we derived analytical expressions for the scaled relative graphs of semimonotone and $\theta$-angle-bounded operators, establishing a connection between the two classes.
We showed that these classes capture the incremental behavior of an Ebers--Moll transistor, enabling us to efficiently compute the response of a nonsmooth and multi-valued common-emitter amplifier circuit with transistors and tunnel diodes using the Chambolle--Pock algorithm.

Future research directions include exploring additional nonmonotone elements like nonlinear capacitors, inductors, and memristors, as well as using the SRG to design circuits with specific input-output behavior.

{\small
\bibliographystyle{elsarticle-num} 
\bibliography{references.bib}
}

\end{document}